\documentstyle[12pt,epsfig,amsmath,amssymb]{article}
\begin{document}

\noindent Stockholm, USITP 04\\
February 2004\\

%\vspace{2mm}

\begin{center}

{\large BIRKHOFF'S POLYTOPE AND}

\vspace{3mm}

{\large UNISTOCHASTIC MATRICES}

\vspace{3mm}

{\large N = 3 AND N = 4}

\vspace{10mm}

Ingemar Bengtsson*\footnote{ingemar@physto.se} \hspace{10mm} \AA sa Ericsson*\footnote{asae@physto.se}
\hspace{10mm} Marek Ku\'s**\footnote{marek@cft.edu.pl} 

\

\hspace{4mm} 
Wojciech Tadej***\footnote{wtadej@wp.pl}
\hspace{15mm} Karol \.Zyczkowski****\footnote{karol@tatry.if.uj.edu.pl}

\vspace{7mm}

*{\it Stockholm University, AlbaNova, Fysikum, 106 91 Stockholm, Sweden.}

**{\it Centrum Fizyki Teoretycznej, Polska Akademia Nauk,
Al. Lotnik\'{o}w 32/44, 02-668 Warszawa, Poland.}

***{\it Cardinal Stefan Wyszynski University, Warszawa, Poland.}

****{\it Instytut Fizyki im. Smoluchowskiego,
Uniwersytet Jagiello\'{n}ski, ul. Reymonta 4, 30-059 Krak\'{o}w, Poland.}

\vspace{8mm}

{\bf Abstract}

\end{center}

\vspace{2mm}

\noindent The set of bistochastic or doubly stochastic $N \times N$ matrices 
form a convex set called Birkhoff's polytope, that we describe in some detail.
Our problem is to characterize the set of unistochastic
matrices as a subset of Birkhoff's polytope. For $N = 3$ we present
fairly complete results. For $N = 4$ partial results are obtained. An
interesting difference between the two cases is that there is a ball
of unistochastic matrices around the van der Waerden matrix for $N = 3$,
while this is not the case for $N = 4$.

\newpage

{\bf 1. Introduction}

\vspace{5mm}

\noindent There is a surprising variety of contexts in which unistochastic
matrices arise, and any one of them may be taken as a motivation for the
present study. But let us first define our terms: An $N \times N$ matrix
$B$ is said to be bistochastic if its matrix elements obey

\begin{equation} \mbox{i:} \ B_{ij} \geq 0 \hspace{10mm} \mbox{ii:} \
\sum_iB_{ij} = 1 \hspace{10mm} \mbox{iii:} \ \sum_jB_{ij} = 1 \ .
\label{1} \end{equation}

\noindent The name ``bistochastic'' has to do with the fact
that these matrices are usually supposed to act on probability distributions,
thought of as $N$ component vectors. The first condition ensures that
positive vectors are transformed to positive vectors, the second that
the sum of the components of the vectors remains invariant, and the third
that the uniform distribution (a vector with all components equal) is
a fixed point of the map. Hence a bistochastic matrix causes a kind of
contraction of the probability simplex towards the uniform distribution.
Condition iii is important because it guarantees that the map increases
entropy. We obtain a bistochastic matrix if we start with a unitary matrix
$U$ and take the absolute value squared of its matrix elements,

\begin{equation} B_{ij} = |U_{ij}|^2 \ . \label{2} \end{equation}

\noindent For connaisseurs of linear algebra, $B$ is the Hadamard product
of $U$ and its complex conjugate. If there exists such a $U$ then $B$
is said to be unistochastic. If $U$ is also real, that is orthogonal,
then $B$ is said to be orthostochastic.
Bistochastic matrices arise frequently in situations where probability
distributions are changing, and we will soon see why one may want them
to be unistochastic.

A somewhat distinguished bistochastic matrix is the van der Waerden matrix
$B_{\star}$, whose matrix elements obey

\begin{equation} B_{{\star}ij} = \frac{1}{N} \ . \end{equation}

\noindent The van der Waerden matrix is unistochastic. A corresponding
unitary matrix is known as a complex Hadamard matrix. An example of a complex
Hadamard matrix is the Fourier matrix, whose matrix elements are

\begin{equation} U_{ij} = \frac{1}{\sqrt{N}}q^{ij} \ , \hspace{8mm} 0
\leq i,j \leq N-1 \ . \end{equation}

\noindent Here $q = e^{2{\pi}i/N}$ is a root of unity. Complex Hadamard
matrices have a long history in mathematics \cite{Sylvester}. If the matrix
$U$ is also real it is referred to simply as a Hadamard matrix. By the way,
the name of the van der Waerden matrix has to do with the conjecture that
this matrix has the largest permanent of all bistochastic matrices;
this is true but took a long time to prove \cite{vanLint}.

It is clear that we now have two mathematical questions on our hands:

{\bf I}: Given a bistochastic matrix, is it unistochastic?

{\bf II}: If so, to what extent is $U$ determined by $B$?

\noindent The answer to question ${\bf I}$ turns out to depend on
the bistochastic matrix chosen, so that in effect question ${\bf I}$ turns
into the problem of characterizing the unistochastic subset of the
set of all bistochastic matrices. But why are these questions
interesting? The answer is that they naturally turn up in many contexts.
Let us give a partial list of those.

The first context has to do with the foundations of quantum mechanics.
Here there are a number of approaches where one begins by arguing that
transition probabilities, suitably defined, form bistochastic matrices.
In attempting to build some group structure into these transition
probabilities one is then led to require that they form unistochastic
matrices, and so one runs into question ${\bf I}$. A sample of
the literature includes Land\'e \cite{Lande}, Rovelli \cite{Rovelli} and
Khrennikov \cite{Khrennikov}.

The second context is classical computer science, especially the theory of
error correcting codes, design theory, and other areas of discrete
mathematics where real Hadamard matrices are very useful.
Hadamard conjectured that such matrices exist when $N = 2$ and
$N = 4k$ \cite{Hadamard}. The conjecture is still open, although
much is known \cite{Hedayat}. For explicit examples of real Hadamard
matrices of all orders up to $256 \times 256$, consult Sloane \cite{Sloane}.

The third context is quantum information theory, where the restriction to
real Hadamard matrices is somewhat unnatural. Complex Hadamard
matrices have been studied in the quantum optics community in
the guise of symmetric multiports; they are examples of specially designed
unitary transformations that can be realized in the laboratory
\cite{Zeilinger} \cite{Stenholm}. There is also an interesting theorem
\cite{Werner} to the effect that the classification of all possible
teleportation schemes can be reduced to the classification of all
sets of maximally entangled bases; this is relevant here because
such sets can be obtained from the combination of a Latin square
and a complex Hadamard matrix. The construction of all possible
Latin squares has nothing to do with us here, but the construction
of all complex Hadamard matrices certainly has. Mathematicians
have studied this problem with various motivations \cite{Bjorck}
\cite{Haagerup} \cite{Petrescu} \cite{Dita}.

The fourth context is the attempt to formulate quantum mechanics on
graphs (in the laboratory on thin strips of, say, gold film). Here question
{\bf I} arises as a question about what Markov processes that have a
quantum counterpart in the given setting \cite{Tanner} \cite{Pakonski1}
\cite{Pakonski2}. In this
connection studies of the spectra and entropies of unistochastic matrices
chosen at random have been made \cite{Karol}; we will return to some
of these issues below.

A fifth context is particle physics, where the interest centers on
question {\bf II}. Thus in the theory of weak interactions we encounter
the unitary Kobayashi-Maskawa matrices (one for quarks and one for
neutrinos), and Jarlskog raised the question to what extent such a matrix
can be parametrized by the easily measured moduli of its matrix elements.
The physically interesting case here is $N = 3$ \cite{Jarlskog}, and
possibly also $N = 4$ in case a fourth generation of quarks should be
discovered \cite{Auberson}. (The same question occurs in scattering theory,
with no restriction on $N$ \cite{Mennessier}.)

We end the list of possible applications here, and turn to the organisation
of our paper. In section 2 we describe the set of all bistochastic
matrices ${\cal B}_N$. It is a convex polytope well known to practioneers
of linear programming; it is sometimes called the assignment polytope
because it arises in the problem of assigning $N$ workers to $N$ tasks,
given their efficiency ratings for each task. We describe the cases $N = 3$
and $N = 4$ in detail ($N = 2$ is trivial). In section 3 we discuss some
generalities concerning unistochastic matrices, and then we characterize the
unistochastic subset for the case of $N = 3$. Most of our results can
be found elsewhere but, we believe, not in this coherent form. In
section 4 we address the same question for $N = 4$. This turns out to be
a more difficult task, but at the end of this section there will be a proof
that every neighbourhood of the van der Waerden matrix contains matrices
that are not unistochastic. This is a striking difference to the $N = 3$ case.
Section 5 summarises our conclusions. Some technical matters are found in
three appendices. Our results on $N > 4$ will be reported in a separate
publication \cite{Tadej}.

\newpage

{\bf 2. Birkhoff's polytope}

\vspace{5mm}

\noindent The set ${\cal B}_N$ of bistochastic $N \times N$ matrices
has $(N-1)^2$ dimensions. To do this count, note that the last row
and the last column are fixed by the conditions that the row and
column sums should equal one. The remaining $(N-1)^2$ entries can be chosen
freely, within limits. Birkhoff proved that ${\cal B}_N$ is a convex
polytope whose extreme points, or corners, are the $N!$ permutation
matrices \cite{Birkhoff}. All corners are equivalent in the
sense that they can be taken into each other by means of orthogonal
transformations. A bistochastic matrix belongs to the boundary of ${\cal B}_N$
if and only if one of its entries is zero. The boundary consists of
corners, edges, faces, 3-faces and so on; the highest dimensional
faces are called facets and consist of matrices with only one zero
entry. For a detailed account of ${\cal B}_N$, especially its
face structure, see Brualdi et al. \cite{Brualdi}. We will be even
more detailed concerning ${\cal B}_3$ and ${\cal B}_4$. For definiteness
all 24 permutation matrices that occur when $N = 4$ are listed in
Appendix A.

It is convenient to regard the convex polytope ${\cal B}_N$ as a subset
of a vector space with the van der Waerden matrix $B_{\star}$ as its
origin. The distance squared between two matrices is chosen to be

\begin{equation}
D^2(A,B) = \mbox{Tr}(A-B)(A^{\dagger} - B^{\dagger}) \ .
\label{dist}
 \end{equation}

\noindent The distance squared between an arbitrary bistochastic matrix
$B$ and $B_{\star}$ is then given by

\begin{equation} D^2(B,B_{\star}) = \sum_{i,j} B^2_{ij} - 1 \ . \end{equation}

\noindent In particular, the distance between $B_{\star}$ and a corner of
the polytope becomes $D = \sqrt{N-1}$. Permutations of rows or columns
are orthogonal transformations since they preserve distance and leave
$B_{\star}$ invariant. They also take permutation matrices (corners)
into permutation matrices, hence they are symmetry operations of Birkhoff's
polytope as well.

The (Shannon) entropy of a
bistochastic matrix is defined as the entropy of the rows averaged over
the columns,

\begin{equation} S = - \frac{1}{N}\sum_i\sum_jB_{ij}\ln{B_{ij}} \ . \end{equation}

\noindent Its maximum value $\ln{N}$ is attained at $B_{\star}$.
For some of its properties consult S{\l}omczy\'nski \cite{WS} et al. \cite{Karol}.

When $N = 2$ there are just two permutation matrices and ${\cal B}_2$ is
a line segment between these two points. A general bistochastic matrix
can be parametrized as

\begin{equation} B = \left[ \begin{array}{cc} c^2 & s^2 \\ s^2 & c^2 \end{array}
\right] \ , \hspace{8mm} c \equiv \cos{\theta} \ , \hspace{5mm} s \equiv
\sin{\theta} \ , \hspace{3mm} 0 \leq {\theta} \leq \frac{\pi}{2} \ .
\label{8} \end{equation}

\noindent When $N =3$ we have six permutation
matrices forming the vertices of
a four dimensional polytope. It admits a simple description:

\

\noindent \underline{Theorem 1}: The 6 corners of ${\cal B}_3$ are the corners
of two equilateral triangles placed in two orthogonal 2-planes and centered
at $B_{\star}$.

\

\noindent The proof is easy, using as corners the permutation matrices $P_0, P_1,
\dots, P_5$ from Appendix A. The two equilateral triangles are the convex
combinations

\begin{equation} {\Delta}_1 = p_0P_0 + p_3P_3 + p_4P_4 =
\left[ \begin{array}{ccc} p_0 & p_3 & p_4 \\
p_4 & p_0 & p_3 \\ p_3 & p_4 & p_0 \end{array} \right] \ ,
\hspace{5mm} p_0 + p_3 + p_4 = 1 \hspace{2mm} \end{equation}

\noindent and

\begin{equation} {\Delta}_2 = p_1P_1 + p_2P_2 + p_5P_5 =
\left[ \begin{array}{ccc} p_1 & p_2 & p_5 \\
p_2 & p_5 & p_1 \\ p_5 & p_1 & p_2 \end{array} \right] \ ,
\hspace{5mm} p_1 + p_2 + p_5 = 1 \ . \end{equation}

\noindent The calculation we have to do is to check that $D^2(P_0, P_3) =
D^2(P_0,P_4) = D^2(P_3, P_4) = 6$ and similarly for the other triangle,
and also that

\begin{equation} \mbox{Tr}({\Delta}_1 - B_{\star})({\Delta}_2^{\dagger} -
B_{\star}^{\dagger}) = 0\end{equation}

\noindent for all values of $p_i$. This is so.

\begin{figure}
        \centerline{ \hbox{
                \epsfig{figure=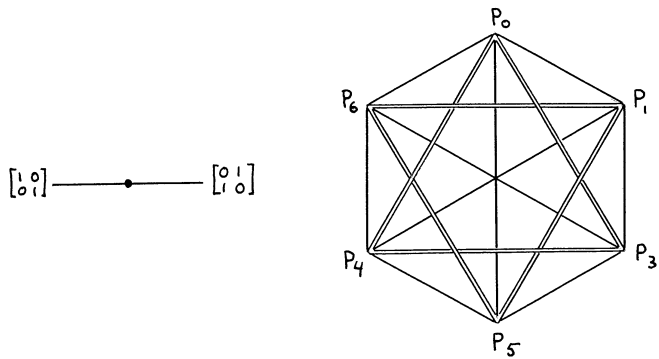,width=9cm}}}
        \caption{Left: Birkhoff's polytope for $N = 2$ (centered at $B_{\star}$).
Right: The graph of Birkhoff's polytope for $N = 3$; single lines have $D^2 = 4$
and double $D^2 = 6$. The double edges form the triangles mentioned in Theorem
1.}
        \label{fig:1}
\end{figure}

\begin{figure}
        \centerline{ \hbox{
                \epsfig{figure=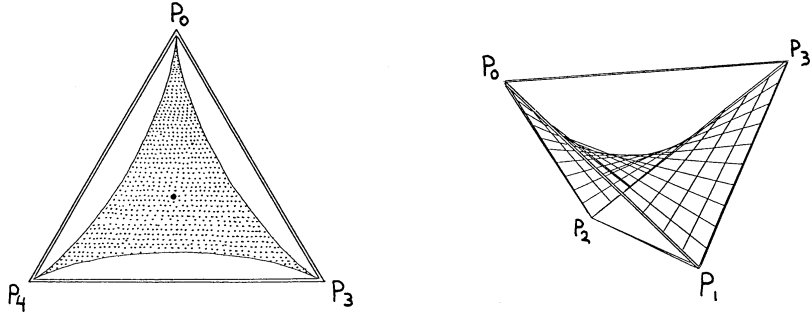,width=11cm}}}
        \caption{Birkhoff's polytope for $N = 3$. Left: One of the two orthogonal
equilateral triangles centered at $B_{\star}$, with its unistochastic subset (the
boundary is the famous hypocycloid). Right: A facet, an irregular tetrahedron,
with its doubly ruled surface of unistochastic matrices.}
        \label{fig:2}
\end{figure}

There are thus 6 corners and $6\cdot 5/2 = 15$ edges, all of which are
extremal. The last is a rather exceptional property; in 3 dimensions
only the simplex has it. There are 9 short edges of length squared $D^2 = 4$
and 6 long edges of length squared $D^2 = 6$, namely the sides of the two
equilateral triangles. A useful overview of ${\cal B}_3$ is given by its
graph, where we exhibit all corners and all edges (see fig. \ref{fig:1}). All
the 2-faces are triangles with one long and two short edges. The 3-faces
in a 4 dimensional polytope are called facets and are made of matrices with
a single zero. They are irregular tetrahedra with two long edges, one
from each equilateral triangle (see fig. \ref{fig:2}).

The volume of ${\cal B}_3$ is readily computed because it can be
triangulated using only three simplices. The total volume is $9/8$.
As $N$ grows the total volume of ${\cal B}_N$ becomes increasingly
hard to compute; mathematicians know it for $N \leq 10$ \cite{Beck}.

The next case is the 9 dimensional polytope ${\cal B}_4$. It has 24 corners
and 276 edges. The latter come in four types and we give the classification
including the angle they subtend at ${\cal B}_{\star}$ and whether they
consist of unistochastic matrices or not (see sections 3 and 4):

\begin{equation} \begin{array}{lclll} \ & \mbox{Length squared} &
\mbox{Unistochastic} & \mbox{Angle at origin} & \mbox{Number of edges} \\
4U & 4 & \mbox{Yes} & \mbox{Acute} & 72 \\ 6 & 6 & \mbox{No} & 90 \
\mbox{degrees} & 96 \\ 8 & 8 & \mbox{No} & \mbox{Obtuse} & 72 \\
8U & 8 & \mbox{Yes} & \mbox{Obtuse} & 36 \end{array} \nonumber
\label{12} \end{equation}

\noindent All edges except the $8U$ ones are extremal. The 2-faces consist
of triangles and squares. (Interestingly, for all $N$ it is true that
the 2-faces of Birkhoff's polytope ${\cal B}_N$ are either triangles or
rectangles \cite{Brualdi}.) There are 18 squares
bounded by edges of type $4U$ and their diagonals are of type $8U$. Three
squares meet at each corner. If we pick four permutation matrices we obtain
a 3-face, with six exceptions. The exceptions form 6 regular tetrahedra
centered at $B_{\star}$ whose edges are non-extremal $8U$ edges. They are
denoted $T_i$ and explicitly listed in Appendix A; an example is

\begin{equation} T_1 = p_0P_0 + p_7P_7 + p_{16}P_{16} + p_{23}P_{23} =
\left[ \begin{array}{cccc} p_0 & p_7 & p_{16} & p_{23} \\
p_7 & p_0 & p_{23} & p_{16} \\ p_{16} & p_{23} & p_0 & p_7 \\
p_{23} & p_{16} & p_7 & p_0 \end{array} \right] \ . \label{13} \end{equation}

\begin{figure}
        \centerline{ \hbox{
                \epsfig{figure=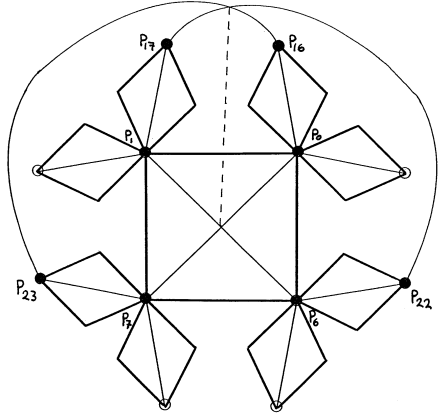,width=6cm}}}
        \caption{How to begin to draw the surface of ${\cal B}_4$.
Two tetrahedra whose edges are the non-extremal diagonals of squares
are shown. The dashed line goes through the polytope; it connects the midpoints
of two opposing $8U$ edges of two tetrahedra that are otherwise disjoint.}
        \label{fig:3}
\end{figure}

\noindent When regular tetrahedra are mentioned below it is understood
that we refer to one of these six. In a sense the structure can now be
drawn; see fig. \ref{fig:3}. The
facets consist of matrices with one zero, so there are 16 facets.

A subset of ${\cal B}_4$ that has no counterpart for ${\cal B}_3$ is the
set of matrices that are tensor products of two by two bistochastic matrices.
This subset splits naturally into several two dimensional components,
and it turns out that they sit in ${\cal B}_4$ as doubly ruled surfaces
inside the regular tetrahedra. Thus the following matrix, parametrised with
two angles, is a tensor product of two matrices of the form (\ref{8}):

\begin{equation} \left[ \begin{array}{cccc} c_1^2c_2^2 & c_1^2s_2^2 &
s_1^2c_2^2 & s_1^2s_2^2 \\ c_1^2s_2^2 & c_1^2c_2^2 & s_1^2s_2^2 & s_1^2c_2^2 \\
s_1^2c_2^2 & s_1^2s_2^2 & c_1^2c_2^2 & c_1^2s_2^2 \\
s_1^2s_2^2 & s_1^2c_2^2 & c_1^2s_2^2 & c_1^2c_2^2 \end{array} \right] \ ,
\hspace{8mm} c_1 \equiv \cos{{\theta}_1} \ \mbox{etc.}\end{equation}

\noindent These matrices form a doubly ruled surface inside the regular
tetrahedron (\ref{13}), analogous to that depicted in fig. \ref{fig:2}.

\begin{figure}
        \centerline{ \hbox{
                \epsfig{figure=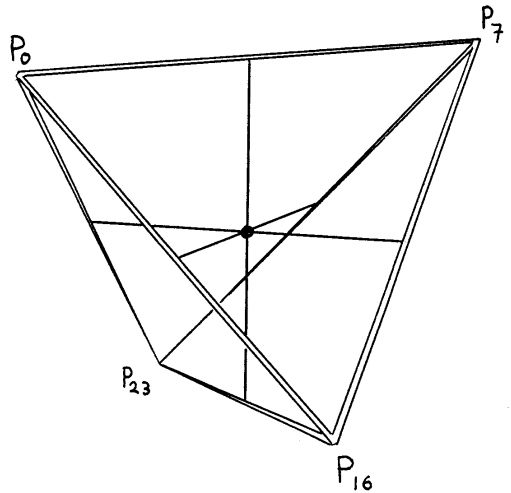,width=7cm}}}
        \caption{A regular tetrahedron centered at $B_{\star}$. It contains the
normal vectors of three orthogonal hyperplanes and belongs entirely to
another six. There are six such regular tetrahedra and pairs of them intersect along
the normal vectors they contain. (Note that the dashed line in Fig. \ref{fig:3}
represents such a normal vector.)}
        \label{fig:5}
\end{figure}

An interesting way to view ${\cal B}_4$, and one that will recur in section
4, stems from the following observation:

\

\noindent \underline{Theorem 2:} The 24 corners of ${\cal B}_4$ belong to a
set of nine orthogonal hyperplanes through ${\cal B}_{\star}$. Each
regular tetrahedron belongs to six hyperplanes and contains the normal
vectors of the remaining three hyperplanes. Each hyperplane contains four
regular tetrahedra and its normal vector is the intersection of the
remaining two regular tetrahedra.

\

\noindent Again the proof is a simple calculation, once the explicit
form of the hyperplanes is known. They are denoted ${\Pi}_i$ and listed
in Appendix A. From now
on, ``hyperplane'' always refers to one of these nine. Fig.
\ref{fig:5} in a sense illustrates the theorem.

It is quite helpful to have an incidence table for tetrahedra and hyperplanes
available. It is

\begin{equation} \begin{array}{cccccccccc} \ & {\Pi}_1 & {\Pi}_2 & {\Pi}_3 &
{\Pi}_4 & {\Pi}_5 & {\Pi}_6 & {\Pi}_7 & {\Pi}_8 & {\Pi}_9 \\
T_1 & \ & X & X & X & \ & X & X & X & \ \\
T_2 & \ & X & X & X & X & \ & X & \ & X \\
T_3 & X & \ & X & \ & X & X & X & X & \ \\
T_4 & X & X & \ & \ & X & X & X & \ & X \\
T_5 & X & \ & X & X & X & \ & \ & X & X \\
T_6 & X & X & \ & X & \ & X & \ & X & X \\
\end{array} \label{15} \end{equation}

\noindent where the tetrahedra $T_i$ and the hyperplanes ${\Pi}_i$ are listed
in Appendix A.

For later purposes we will need some information about exactly how the
hyperplanes divide the space into $2^9$ hyperoctants. For this reason
we look at the rays

\begin{equation} B_i(t) = B_{\star} + tV_i \ , \label{ray} \end{equation}

\noindent where $V_i$ is a vector constructed in terms of the normal vectors
$n_1, \dots, n_9$ of the hyperplanes (see Appendix B), namely

\begin{equation} V_1 \equiv n_1 + n_2 + n_3 + n_4 + n_5 + n_6 + n_7 + n_8 + n_9 =
\frac{1}{4}\left[ \begin{array}{rrrr} 9 & - 3 & - 3 & - 3 \\ - 3 & 1 & 1 & 1 \\
- 3 & 1 & 1 & 1 \\ - 3 & 1 & 1 & 1 \end{array} \right] \label{17} \end{equation}

\begin{equation} V_2 \equiv n_1 + n_2 + n_3 + n_4 + n_5 + n_6 + n_7 + n_8 - n_9 =
\frac{1}{4}\left[ \begin{array}{rrrr} 7 & - 1 & - 1 & - 5 \\ - 1 & - 1 & - 1 & 3 \\
- 1 & - 1 & - 1 & 3 \\ - 5 & 3 & 3 & - 1 \end{array} \right] \end{equation}

\begin{equation} V_3 \equiv n_1 + n_2 + n_3 + n_4 - n_5 + n_6 + n_7 + n_8 - n_9 =
\frac{1}{4}\left[ \begin{array}{rrrr} 5 & 1 & - 3 & - 3 \\ 1 & - 3 & 1 & 1 \\
- 3 & 1 & -3 & 5 \\ - 3 & 1 & 5 & -3 \end{array} \right] \ . \label{19} \end{equation}

\noindent All other
cases can be obtained from one of these three by permutations of rows and
columns. The various hyperoctants are convex cones centered on these rays.
This gives a classification of the hyperoctants into six different types
(since the parameter $t$ can be positive or negative) called respectively
type $\mbox{I}_{\pm}$, $\mbox{II}_{\pm}$ and $\mbox{III}_{\pm}$. Type I
has 16 representatives and is especially noteworthy. For type $\mbox{I}_-$ the
centered ray hits the boundary in the center of one of the 16 facets, at
the matrix $B_1(- \frac{1}{9})$. In the other direction we also hit quite
distinguished points. There are 16 ways of setting one entry of a bistochastic
matrix equal to one, and this gives rise to 16 copies of ${\cal B}_3$ sitting
in the boundary of ${\cal B}_4$. For the octants $\mbox{I}_+$ the centered
ray hits the boundary precisely at the center of such a ${\cal B}_3$, at
the matrix $B_1(\frac{1}{3})$.

In section 4 we will see how
the structure of the unistochastic subset is related to the structure
of Birkhoff's polytope, and in particular to the features we have
stressed.

\newpage

{\bf 3. The unistochastic subset, mostly $N = 3$}

\vspace{5mm}

\noindent Let us begin with some generalities concerning the unistochastic
subset ${\cal U}_N$ of ${\cal B}_N$. The dimension of ${\cal B}_N$ is $(N-1)^2$
and the dimension of $U(N)$ is $N^2$. Therefore the map $U(N) \rightarrow {\cal B}_N$
cannot be one-to-one. Now it is clear that multiplying a row or a column
by a phase factor---an operation that we refer to as rephasing---will result
in the same bistochastic matrix via eq.
(\ref{2}). Therefore the map is naturally defined as a map from a double
coset space to ${\cal B}_N$. The double coset space is

\begin{equation} U(1)\times \dots \times U(1)\setminus U(N)\ / \ U(1)\times \dots \times
U(1) \ , \end{equation}

\noindent with $N$ $U(1)$ factors on the right and $N -1$ factors on the left,
say. The dimension of this set is $(N-1)^2$ so now the dimensions match.
There is a complication because the double coset space is not a smooth manifold.
The action from the left of the $U(1)$ factors on the right coset space
(in itself a well behaved flag manifold) has fixed points. These fixed
points are easy to locate however (and always map to the boundary of
${\cal B}_N$), so that for most practical purposes
we can think of our map as a map between smooth manifolds.

In general we will see that the image of our map is a proper subset of
${\cal B}_N$, and the map is many-to-one. There is not
much we can usefully say about the general case, except for two remarks:
The unistochastic subset ${\cal U}_N$ has the full dimension $(N-1)^2$
while the unistochastic subset of the boundary of ${\cal B}_N$ has
dimension $(N-1)^2 - 2$; why this is so will presently become clear.

For $N = 2$ every bistochastic matrix is orthostochastic. A unitary
matrix that maps to the matrix in eq. (\ref{8}) is

\begin{equation} U = \left[ \begin{array}{cc} c & s \\ s & - c \end{array}
\right] \ , \hspace{8mm} c \equiv \cos{\theta} \ , \hspace{5mm} s \equiv
\sin{\theta} \ , \hspace{3mm} 0 \leq {\theta} \leq \frac{\pi}{2} \ . \end{equation}

\noindent The matrix is given in dephased form. This means that the
first row and the first column is real and positive. This fixes the $U(1)$
factors mentioned above (unless there is a zero entry in one of these
places) and from now on we will present all unitary
matrices in this form. For any $N$ it is straightforward to check whether a given
edge of ${\cal B}_N$ is unistochastic. For $N = 3$ the edges of length
squared equal to 4 are unistochastic, and for $N = 4$ we have the
results given in table (\ref{12}).

Given a $3 \times 3$ bistochastic matrix it is easy to check whether
it is unistochastic or not \cite{Poon} \cite{Jarlskog}. We form the
moduli $r_{ij} = \sqrt{B_{ij}}$ and write down the matrix

\begin{equation} U = \left[ \begin{array}{lll} r_{00} & r_{01} & \bullet \\
r_{10} & r_{11}e^{i{\phi}_{11}} & \bullet \\ r_{20} & r_{21}e^{i{\phi}_{21}}
& \bullet \end{array} \right] \end{equation}

\noindent If this matrix is unitary the original matrix is unistochastic.
The unitarity conditions simply say that the first two columns are orthogonal;
the last column by construction has the right moduli and does not impose
any further restrictions. Therefore the problem is to form a triangle from
three line segments of given lengths

\begin{equation} L_0 = r_{00}r_{01} \hspace{8mm} L_1 = r_{10}r_{11}
\hspace{8mm} L_2 = r_{20}r_{21} \ . \label{23} \end{equation}

\noindent This is possible if and only if the ``chain--links'' conditions
are fulfilled, i.e.

\begin{equation} |L_1 - L_2| \leq L_0 \leq L_1 + L_2 \ . \label{24} \end{equation}

\noindent The bistochastic matrix corresponding to $U$ sits at the
boundary of ${\cal U}_3$ if and only if one of these inequalities is 
saturated. When eq. (\ref{24}) holds the solution is

\begin{equation} \cos{{\phi}_{11}} = \frac{L_2^2 - L_0^2 - L_1^2}{2L_0L_1}
\hspace{6mm} \cos{{\phi}_{21}} = \frac{L_1^2 - L_2^2 - L_0^2}{2L_0L_2}
\end{equation}

\begin{equation} \cos{({\phi}_{11} - {\phi}_{21})} = \frac{L_0^2 -
L_1^2 - L_2^2}{2L_1L_2} \ . \end{equation}

\noindent There is a two-fold ambiguity (corresponding to taking the
complex conjugate of the matrix, $U \rightarrow U^*$). The area $A$ of the
triangle is %given by
%
%\begin{equation} 16A^2 = 2L_0^2L_1^2 + 2L_0^2L_2^2 + 2L_1^2L_2^2 -
%L_0^4 - L_1^4 - L_2^4 \ . \end{equation}
%
%\noindent
easily computed and the chain--links conditions are equivalent to the single inequality
$A \geq 0$. As a matter of fact we can form six so called unitarity triangles
in this way, depending on what pair of columns or rows that we choose.
Although their shapes differ their area is the same, by unitarity \cite{Jarlskog}.

Because we can easily decide if a given matrix is unistochastic, it is easy 
to characterize the unistochastic set ${\cal U}_3$. We single out the following 
facts (some of which are known \cite{Poon}) for attention:

\

\noindent \underline{Theorem 3}: The unistochastic subset ${\cal U}_3$ of
${\cal B}_3$ is a non-convex star shaped four dimensional set whose boundary
consists of the set of orthostochastic matrices. It contains a unistochastic
ball of maximal radius $\sqrt{2}/3$, centered at $B_{\star}$. The set meets
the boundary of ${\cal B}_3$ in a doubly ruled surface in each facet.

\

\noindent The relative
volume of the unistochastic subset is, according to our numerics,

\begin{equation} \frac{\mbox{vol}({\cal U}_3)}
{\mbox{vol}({\cal B}_3)} \approx 0.7520 \pm 0.0005 \ . \end{equation}

\noindent We did not attempt an analytical calculation; details of our
numerics are in Appendix B.

Theorem 3 is easy to prove. To see that ${\cal U}_3$ is
non-convex we just draw its intersection with one of the equilateral
triangles that went into the definition of the polytope, and look at it
(see fig. \ref{fig:2}). An amusing side remark is that the boundary of
the unistochastic set in this picture is a 3-hypocycloid \cite{Karol}.
It can be obtained by rolling a circle of radius $1/3$ inside the unit
circle. The maximal unistochastic ball is centered at $B_{\star}$ and touches
the boundary at the hypocycloid, as one might guess from the picture; its
radius was deduced from results presented in ref. \cite{India}. To see
that the boundary consists of orthostochastic matrices is the observation
that when the chain--links conditions are saturated the phases in $U$
will equal $\pm 1$. That the set is star shaped then follows from an
explicit check that there is only one orthostochastic matrix on any
ray from $B_{\star}$. Finally fig. \ref{fig:2}  includes an explicit
picture of the unistochastic subset of a facet. The reason why it has
codimension one is that a matrix on the boundary of ${\cal B}_N$ has a
zero entry, which means that the number of phases available in the
dephased unitary matrix drops with one, and then the dimension of the
unistochastic set also drops with one; the argument goes through for
any $N$.

Finally let us make some remarks on entropy. We compare the
Shannon entropy averaged over ${\cal B}_3$ using the flat measure, the
Shannon entropy averaged over ${\cal U}_3$ also using the flat measure,
and the maximal Shannon entropy $S_{\rm max}$. Numerically we find that

\begin{equation}
\langle S\rangle_{{\cal B}_3} \approx 0.883
{\rm \quad and \quad}
\langle S\rangle_{{\cal U}_3} \approx 0.908
\label{entr3}
 \end{equation}

\noindent with all digits significant. Observe that the latter average is
larger since some matrices of small entropy close to the boundary of
${\cal B}_3$ are not unistochastic and do not contribute to the average
over ${\cal U}_3$. The above data may be compared with the maximal possible
entropy $S_{\rm max} = \ln{3} \approx 1.099$, attained at $B_{\star}$, and
also with

\begin{equation} \langle S\rangle_{\rm Haar} = \frac{1}{2} + \frac{1}{3}
\approx 0.833 \ , \end{equation}

\noindent which is the average taken over ${\cal U}_3$ with respect to
measure induced by the Haar measure on $U(3)$. This analytical result follows
directly from the work of Jones, who computed the average entropy of squared
components of complex random vectors \cite{Jones}. It is easy to see that
the two averages coincide. For details of our numerics
consult Appendix B.

\newpage

{\bf 4. The unistochastic subset, mostly $N = 4$}

\vspace{5mm}

\noindent The case $N = 4$ is more difficult. It is also clear from the
outset that it will be qualitatively different---thus the dimension of
the orthogonal group is too small for the boundary of the unistochastic
set ${\cal U}_4$ to be formed by orthostochastic matrices alone. There
are other differences too, as we will see.

Given a bistochastic matrix we can again define $r_{ij} = \sqrt{B_{ij}}$
and consider

\begin{equation} U = \left[ \begin{array}{llll} r_{00} & r_{01} & r_{02} &\bullet \\
r_{10} & r_{11}e^{i{\phi}_{11}} & r_{12}e^{i{\phi}_{12}} & \bullet \\
r_{20} & r_{21}e^{i{\phi}_{21}} & r_{22}e^{i{\phi}_{22}} & \bullet \\
r_{30} & r_{31}e^{i{\phi}_{31}} & r_{32}e^{i{\phi}_{32}} & \bullet
\end{array} \right] \end{equation}

\noindent Phases must now be chosen so that this matrix is unitary,
and more especially so that the three columns we focus on are orthogonal.
Geometrically this is the problem of forming three quadrilaterals
with their sides given and six free angles. This is not a simple problem,
and in practice we have to resort to numerics to see whether a given
bistochastic matrix is unistochastic (see Appendix B for details).
There are some easy special cases though. One easy case is that of
a matrix belonging to the boundary of ${\cal B}_N$. Then the matrix
$U$ must contain one zero entry and when we check the orthogonality
of our three columns two of the equations reduce to the problem of
forming triangles, so that the angles are completely fixed when we
consider the final orthogonality relation. Another easy case concerns
the regular tetrahedra. They turn out to consist of orthostochastic
matrices; for the example given in eq. (\ref{13}) a corresponding
orthonormal matrix is

\begin{equation} O_1 = \left[ \begin{array}{rrrr} \sqrt{p_0} & \sqrt{p_7} &
\sqrt{p_{16}} & \sqrt{p_{23}} \\ \sqrt{p_7} & - \sqrt{p_0} &
- \sqrt{p_{23}} & \sqrt{p_{16}} \\ \sqrt{p_{16}} & \sqrt{p_{23}} &
- \sqrt{p_0} & - \sqrt{p_7} \\ \sqrt{p_{23}} & - \sqrt{p_{16}} &
\sqrt{p_7} & - \sqrt{p_0} \end{array} \right] \ . \end{equation}

\noindent This saturates a bound saying that the maximum number of
$N\times N$ permutation matrices whose convex hull is unistochastic is
not larger than $2^{\left[\frac{N}{2}\right]}$, where $\left[ N/2\right]$
denotes the integer part of $N/2$ \cite{Au-Yeung}.

Let us now turn our attention to $B_{\star}$. Hadamard \cite{Hadamard}
observed that up to permutations of rows and columns the most general
form of the complex Hadamard matrix is

\begin{equation} H({\phi}) = \frac{1}{2}\left[ \begin{array}{cccc}
1 & 1 & 1 & 1 \\
1 & e^{i{\phi}} & - 1 & - e^{i{\phi}} \\
1 & - 1 & 1 & - 1 \\
1 & - e^{i{\phi}} & - 1 & e^{i{\phi}}
\end{array} \right] \label{30} \end{equation}

\noindent One can show that this is a geodesic in $U(N)$. The news,
compared to $N = 3$, is that $B_{\star}$ is
orthostochastic because $H(0)$ is real. Moreover there is a continuous
set of dephased unitaries mapping to the same $B$. In a calculational
{\it tour de force}, Auberson et al. \cite{Auberson} were able to
determine all bistochastic matrices whose dephased unitary preimages
contain a continuous ambiguity (and they found that the ambiguity is
given by one parameter in all cases). There
are three such families. Using the notation of ref. \cite{Auberson}
they consist of matrices of the following form:

\begin{equation} \mbox{Type A:} \hspace{5mm} \left[ \begin{array}{cccc}
a & b & c & d \\ b & a & d & c \\ e & f & g & h \\ f & e & h & g
\end{array} \right] \hspace{8mm} \mbox{Type C:} \hspace{5mm} \left[ \begin{array}{cccc}
a & a & \frac{1}{2} - a & \frac{1}{2} - a \\ b & b & \frac{1}{2} - b &
\frac{1}{2} - b \\ c & c & \frac{1}{2} - c & \frac{1}{2} - c \\ d & d &
\frac{1}{2} - d & \frac{1}{2} - d \end{array} \right]
\end{equation}

\begin{equation} \mbox{Type B:} \hspace{5mm}
\left[ \begin{array}{cccc} s_1^2s_2^2 & c_1^2s_2^2 &
c_3^2c_2^2 & s_3^2c_2^2 \\ s_1^2c_2^2 & c_1^2c_2^2 & c_3^2s_2^2 & s_3^2s_2^2 \\
c_1^2c_4^2 & s_1^2c_4^2 & s_3^2s_4^2 & c_3^2s_4^2 \\
c_1^2s_4^2 & s_1^2s_4^2 & s_3^2c_4^2 & c_3^2c_4^2 \end{array} \right] \ .
\hspace{3mm} \end{equation}

\noindent Here $c_1 = \cos{{\theta}_1}$, $s_1 \equiv \sin{{\theta}_1}$,
and so on. Type A consists
of nine five dimensional sets, type B of nine four dimensional sets,
and type C of six three dimensional sets. In trying to understand their
location in ${\cal B}_4$ the observation in section 2 concerning the nine
orthogonal hyperplanes begins to pay dividends. (In particular, consult the
incidence table \ref{15}.) Type A consists of the
linear subspaces obtained by taking all intersections of four hyperplanes
that contain exactly two regular tetrahedra. Type C consists of the linear
subspaces obtained by taking all intersections of six hyperplanes that
contain no permutation matrices at all. Type B finally consists of curved
manifolds confined to one hyperplane. Auberson's families are not exclusive.
In particular tensor product matrices belong to families A and B, which means
that there are two genuinely different ways of introducing a free phase in
the corresponding unitary matrix. Outside the three sets $A$, $B$ and $C$
Auberson et al. find a 12-fold discrete ambiguity in the dephased
unitaries, dropping to 4-fold for symmetric matrices \cite{Auberson}.

Tensor product matrices $B_4 = B_2\otimes B'_2$ appear because $4 = 2\times 2$
is a composite number. That they are always unistochastic follows from a
more general result:

\

\noindent \underline{Lemma 1}: Let $B_K$ and $B_M$ be unistochastic matrices of
size $K$ and $M$, respectively. Then the matrix $B_N = B_K \otimes B_M$ of
size $KM$ is unistochastic. The corresponding
unitary matrices contain at least $(K-1)(M-1)$ free phases when dephased.

\

\noindent That $B_N$ is unistochastic follows from properties of the Hadamard
and the tensor products. By definition, the Hadamard product $A\circ B$ of two
matrices is the matrix whose matrix elements are the products of the corresponding
matrix elements of $A$ and $B$. Then $B_K = U_K\circ U^*_K$ and $B_M =
U_M\circ U^*_M$ implies that $B_N = (U_K\circ U^*_K)\otimes
(U_M\circ U^*_M) = (U_K\otimes U_M)\circ (U^*_K\otimes U^*_M)$, so
it is unistochastic.
The existence of free phases is an easy generalization
of proposition 2.9 in Haagerup \cite{Haagerup}.

The hyperplane structure of ${\cal B}_4$ reverberates in the structure of
the unistochastic set in several ways. Let us consider how the tangent space
of $U(N)$ behaves under the map to ${\cal B}_N$. In equations, this means
that we fix a unitary matrix $U_0$ and expand

\begin{equation} U(t) = e^{iht}U_0 = (1 + iht - \frac{1}{2}h^2t^2 + \dots )
U_0 \end{equation}

\noindent where $h$ is an Hermitian matrix.
Then we study bistochastic matrices with elements $B_{ij}(t) =
|U_{ij}(t)|^2$ to first order in $t$. The following general features are observed:

\begin{itemize}

\item{Generically the tangent space of $U(N)$ maps onto the tangent space of
${\cal B}_N$. This implies that the dimension of the unistochastic set is
equal to that of ${\cal B}_N$; we checked this statement by generating
unitary matrices at random using the Haar measure on the group.}

\item{A matrix element in
$B$ receives a first order contribution only if it is
non-vanishing. Hence the map of the tangent space of $U(N)$ to the tangent
space of ${\cal B}_N$ is degenerate at the boundary of the polytope. In
general such behaviour is to be expected at the boundary of the unistochastic
set ${\cal U}_N$.}

\item{If $U_0$ is real the map is degenerate in the sense that the tangent space
maps to an $N(N-1)/2$ dimensional subspace of the tangent space of ${\cal B}_N$.}

\item{If $U_0$ maps to a corner of the polytope then the first order contributions
vanish. To second order we pick up the tip of a convex cone whose extreme rays
are the $N(N-1)/2$ $4U$ edges emanating from that corner.}

\end{itemize}

For $N = 4$ the story becomes interesting when we choose $U_0$ equal to the Hadamard
matrix $H({\phi})$. Then we find that the tangent space at $U_0$ maps into
one of the nine hyperplanes; which particular one depends on how we permute
rows and columns in eq. (\ref{30}). The question therefore arises whether
the orthostochastic van der Waerden matrix belongs to the boundary of the
unistochastic set---or not since {\it a priori} such degeneracies can occur
also in the interior of the set.

We know that we can form curves of unistochastic matrices starting from $B_{\star}$
and moving out into the nine hyperplanes. Can we form such curves that go directly
out into one of the $2^9$ hyperoctants? Here the division of the $2^9$ hyperoctants
into six different types becomes relevant. We have investigated whether their
central rays given in eqs. (\ref{ray}-\ref{19}) consist of unistochastic matrices,
or not. Let us begin with the 16 hyperoctants of type $\mbox{I}$, where the central ray
$B_1(t) = B_{\star} + tV_1$ hits the boundary in the center of one of the 16
${\cal B}_3$ sitting in the boundary (at $t = 1/3$), and in the center of one
of the 16 facets (at $t = - 1/9$). Of these two points, the first is unistochastic,
the second is not. A one parameter family of candidate unitary matrices
that maps to the central ray is

\begin{equation} U(t) = \frac{1}{2}\left[ \begin{array}{llll} \sqrt{1-3t} &
\sqrt{1-3t} & \sqrt{1-3t} & \bullet \\ \sqrt{1+t} & \sqrt{1+t}e^{{\phi}_{11}}
& \sqrt{1+t}e^{i{\phi}_{12}} & \bullet \\
\sqrt{1+t} & \sqrt{1+t}e^{i{\phi}_{21}} & \sqrt{1+t}e^{{\phi}_{22}} & \bullet \\
\sqrt{1+t} & \sqrt{1+t}e^{{\phi}_{31}} & \sqrt{1+t}e^{i{\phi}_{32}} & \bullet
\end{array} \right] \ , \end{equation}

\noindent where $t>0$ and we permuted the columns relative to eq. (\ref{17})
in order to get the unitarity equations in a pleasant form. (We do not need
to give the phases for the last column.) The conditions
that the first three columns be orthogonal read

\begin{eqnarray} e^{i{\phi}_{11}} + e^{i{\phi}_{21}} + e^{i{\phi}_{31}} + L = 0
\label{37} \\
e^{i{\phi}_{12}} + e^{i{\phi}_{22}} + e^{i{\phi}_{32}} + L = 0 \label{38} \\
e^{i({\phi}_{11}- {\phi}_{12})} + e^{i({\phi}_{21} - {\phi}_{22})} +
e^{i({\phi}_{31} - {\phi}_{32})} + L = 0 \label{39} \end{eqnarray}

\noindent where

\begin{equation} L = \frac{1 - 3t}{1+t} \ . \end{equation}

\noindent In Appendix C we prove that the system of equations (\ref{37}-\ref{39})

\noindent 1. has no real solutions for $L > 1$,

\noindent 2. for $0 < L < 1$ has the solution
\begin{equation} \begin{array}{lll} {\phi}_{11} = 0 & {\phi}_{21} = {\phi} &
{\phi}_{31} = - {\phi} \\ {\phi}_{12} = {\phi} & {\phi}_{22} = 0 & {\phi}_{32} =
- {\phi} \end{array} \ , \hspace{6mm} \cos{\phi} = \frac{t-1}{t+1} =
- \frac{L+1}{2}  \ . \label{41} \end{equation}

\noindent It follows that the central ray is unistochastic for the hyperoctants
of type $\mbox{I}_+$ (and the
unitary matrices on the central ray tend to the real Hadamard matrix at $t = 0$).
In the other direction the central ray is not unistochastic for type $\mbox{I}_-$.
Thus we have proved

\

\noindent \underline{Theorem 4}: For $N = 4$ there are non-unistochastic
matrices in every neighbourhood of the van der Waerden matrix $B_{\star}$. At
$B_{\star}$ the map $U(4) \rightarrow {\cal B}_4$ aligns the tangent space
of $U(4)$ with one of the nine orthogonal hyperplanes.

\

\noindent The structure of the
unistochastic set is dramatically different depending on whether $N = 3$ or
$N = 4$. It is only in the former case that there is a ball of
unistochastic matrices surrounding $B_{\star}$. On the other hand, the
hyperoctants are not empty---some of them do contain unistochastic matrices
all the way down to $B_{\star}$.

Concerning the other hyperoctants, for types $\mbox{II}_-$, $\mbox{III}_+$, and
$\mbox{III}_-$ the central rays hit the boundary of the
polytope in points that are not unistochastic, but numerically we find
that a part of the ray close to $B_{\star}$ is unistochastic. For type
$\mbox{II}_+$ we hit the boundary in a unistochastic point and
numerically we find the entire ray to be unistochastic. There is still
much that we do not know. We do not know if the hyperoctants of type
$\mbox{I}_-$ are entirely free of unistochastic matrices, nor do we know if
${\cal U}_4$ is star shaped, or what its relative volume may be. What is
clear from the results that we do have is that the global structure of
Birkhoff's polytope reverberates in the structure of the unistochastic
subset in an interesting way---it is a little bit like a nine dimensional
snowflake, because the nine hyperplanes in ${\cal B}_4$ can be found
through an analysis of the behaviour of ${\cal U}_4$ in the neighbourhood
of $B_{\star}$.

\newpage

{\bf 5. Conclusions}

\vspace{5mm}

\noindent Our reasons for studying the unistochastic subset of Birkhoff's
polytope have been summarized in the introduction. Because the problem
is a difficult one we concentrated on the cases $N = 3$ and $N = 4$.
Our descriptions of Birkhoff's polytope for these two cases are
given in Theorems 1 and 2, respectively, and a characterization sufficient
for our purposes of the unistochastic set for $N = 3$ was given in Theorem
3. For $N =4$ the dimension of the unistochastic set is again equal to
that of the polytope itself, but its structure differs dramatically from
the $N = 3$ case. In particular Theorem 4 states that for $N = 4$ there
are non-unistochastic matrices in every neighbourhood of the van der
Waerden matrix. Hence there does not exist a unistochastic ball
surrounding the van der Waerden matrix. We observed that the structure
of the unistochastic set at the center
of the polytope reflects the global structure of the latter in an
interesting way.

It is natural to ask to what extent the difference
between the two cases reflects the fact that 3 is prime while 4 is not.
Although this is not the place to discuss the cases $N > 4$, since
some of us intend to do so in a separate publication \cite{Tadej}, let
us mention that the dimension of the unistochastic set is equal to
that of ${\cal B}_N$ for all values of $N$. On the other hand it is
only in the case of $N$ being a prime number that we have been able
to show that there is a unistochastic ball surrounding the van der
Waerden matrix.

\vspace{1cm}

{\bf Acknowledgements:}

\vspace{5mm}

\noindent We thank G\"oran Bj\"orck, Prot Pako{\'n}ski,
 Wojciech S{\l}omczy\'nski, and Gregor Tanner for
discussions, Petre D\u{i}\c{t}a for email correspondence, and
Uffe Haagerup for supplying us with a copy of Petrescu's thesis.
Financial support from the Swedish Research Council VR, and from 
the Polish Ministry of Scientific Research under grant No 
PBZ-MIN-008/P03/2003, is gratefully acknowledged.

\vspace{1cm}

{\bf Appendix A: Notation}

\vspace{5mm}

\noindent An explicit list of permutation matrices for $N = 4$ is

\begin{equation} {\small P_0 = \left[ \begin{array}{cccc} 1 & 0 & 0 & 0 \\ 0 & 1 & 0 & 0 \\
0 & 0 & 1 & 0 \\ 0 & 0 & 0 & 1 \end{array} \right] \hspace{6mm}
P_1 = \left[ \begin{array}{cccc} 1 & 0 & 0 & 0 \\ 0 & 1 & 0 & 0 \\
0 & 0 & 0 & 1 \\ 0 & 0 & 1 & 0 \end{array} \right] \hspace{6mm}
P_2 = \left[ \begin{array}{cccc} 1 & 0 & 0 & 0 \\ 0 & 0 & 1 & 0 \\
0 & 1 & 0 & 0 \\ 0 & 0 & 0 & 1 \end{array} \right] } \end{equation}

\begin{equation} {\small \hspace{6mm} P_3 = \left[ \begin{array}{cccc}
1 & 0 & 0 & 0 \\ 0 & 0 & 1 & 0 \\
0 & 0 & 0 & 1 \\ 0 & 1 & 0 & 0 \end{array} \right] \hspace{6mm}
P_4 = \left[ \begin{array}{cccc} 1 & 0 & 0 & 0 \\ 0 & 0 & 0 & 1 \\
0 & 1 & 0 & 0 \\ 0 & 0 & 1 & 0 \end{array} \right] \hspace{6mm}
P_5 = \left[ \begin{array}{cccc} 1 & 0 & 0 & 0 \\ 0 & 0 & 0 & 1 \\
0 & 0 & 1 & 0 \\ 0 & 1 & 0 & 0 \end{array} \right] } \end{equation}

\begin{equation} P_6 = \left[ \begin{array}{cccc} 0 & 1 & 0 & 0 \\ 1 & 0 & 0 & 0 \\
0 & 0 & 1 & 0 \\ 0 & 0 & 0 & 1 \end{array} \right] \hspace{6mm}
P_7 = \left[ \begin{array}{cccc} 0 & 1 & 0 & 0 \\ 1 & 0 & 0 & 0 \\
0 & 0 & 0 & 1 \\ 0 & 0 & 1 & 0 \end{array} \right] \hspace{6mm}
P_8 = \left[ \begin{array}{cccc} 0 & 1 & 0 & 0 \\ 0 & 0 & 1 & 0 \\
1 & 0 & 0 & 0 \\ 0 & 0 & 0 & 1 \end{array} \right] \end{equation}

\begin{equation} \hspace{6mm} P_9 = \left[ \begin{array}{cccc}
0 & 1 & 0 & 0 \\ 0 & 0 & 1 & 0 \\
0 & 0 & 0 & 1 \\ 1 & 0 & 0 & 0 \end{array} \right] \hspace{6mm}
P_{10} = \left[ \begin{array}{cccc} 0 & 1 & 0 & 0 \\ 0 & 0 & 0 & 1 \\
1 & 0 & 0 & 0 \\ 0 & 0 & 1 & 0 \end{array} \right] \hspace{6mm}
P_{11} = \left[ \begin{array}{cccc} 0 & 1 & 0 & 0 \\ 0 & 0 & 0 & 1 \\
0 & 0 & 1 & 0 \\ 1 & 0 & 0 & 0 \end{array} \right] \end{equation}

\begin{equation} P_{12} = \left[ \begin{array}{cccc} 0 & 0 & 1 & 0 \\ 1 & 0 & 0 & 0 \\
0 & 1 & 0 & 0 \\ 0 & 0 & 0 & 1 \end{array} \right] \hspace{5mm}
P_{13} = \left[ \begin{array}{cccc} 0 & 0 & 1 & 0 \\ 1 & 0 & 0 & 0 \\
0 & 0 & 0 & 1 \\ 0 & 1 & 0 & 0 \end{array} \right] \hspace{5mm}
P_{14} = \left[ \begin{array}{cccc} 0 & 0 & 1 & 0 \\ 0 & 1 & 0 & 0 \\
1 & 0 & 0 & 0 \\ 0 & 0 & 0 & 1 \end{array} \right] \end{equation}

\begin{equation} \hspace{6mm} P_{15} = \left[ \begin{array}{cccc}
0 & 0 & 1 & 0 \\ 0 & 1 & 0 & 0 \\
0 & 0 & 0 & 1 \\ 1 & 0 & 0 & 0 \end{array} \right] \hspace{5mm}
P_{16} = \left[ \begin{array}{cccc} 0 & 0 & 1 & 0 \\ 0 & 0 & 0 & 1 \\
1 & 0 & 0 & 0 \\ 0 & 1 & 0 & 0 \end{array} \right] \hspace{5mm}
P_{17} = \left[ \begin{array}{cccc} 0 & 0 & 1 & 0 \\ 0 & 0 & 0 & 1 \\
0 & 1 & 0 & 0 \\ 1 & 0 & 0 & 0 \end{array} \right] \end{equation}

\begin{equation} P_{18} = \left[ \begin{array}{cccc} 0 & 0 & 0 & 1 \\ 1 & 0 & 0 & 0 \\
0 & 1 & 0 & 0 \\ 0 & 0 & 1 & 0 \end{array} \right] \hspace{5mm}
P_{19} = \left[ \begin{array}{cccc} 0 & 0 & 0 & 1 \\ 1 & 0 & 0 & 0 \\
0 & 0 & 1 & 0 \\ 0 & 1 & 0 & 0 \end{array} \right] \hspace{5mm}
P_{20} = \left[ \begin{array}{cccc} 0 & 0 & 0 & 1 \\ 0 & 1 & 0 & 0 \\
1 & 0 & 0 & 0 \\ 0 & 0 & 1 & 0 \end{array} \right] \end{equation}

\begin{equation} \hspace{6mm} P_{21} = \left[ \begin{array}{cccc}
0 & 0 & 0 & 1 \\ 0 & 1 & 0 & 0 \\
0 & 0 & 1 & 0 \\ 1 & 0 & 0 & 0 \end{array} \right] \hspace{5mm}
P_{22} = \left[ \begin{array}{cccc} 0 & 0 & 0 & 1 \\ 0 & 0 & 1 & 0 \\
1 & 0 & 0 & 0 \\ 0 & 1 & 0 & 0 \end{array} \right] \hspace{5mm}
P_{23} = \left[ \begin{array}{cccc} 0 & 0 & 0 & 1 \\ 0 & 0 & 1 & 0 \\
0 & 1 & 0 & 0 \\ 1 & 0 & 0 & 0 \end{array} \right] \ . \end{equation}

\noindent In Birkhoff's polytope these 24 matrices form the corners
of 6 regular tetrahedra, namely the convex hulls of the sets

\begin{eqnarray} T_1 = \{P_0,P_7,P_{16},P_{23}\} \hspace{7mm} T_2 =
\{P_1,P_6,P_{17},P_{22}\} \hspace{2mm} \nonumber \\
\ \nonumber \\
T_3 = \{P_2,P_{10},P_{13},P_{21}\} \hspace{6mm}
T_4 = \{P_3,P_{11},P_{12},P_{20}\} \\
\ \nonumber \\
T_5 = \{P_{4},P_{8},P_{15},P_{19}\} \hspace{6mm}
T_6 = \{P_5,P_9,P_{14},P_{18}\} \ . \nonumber \end{eqnarray}

\noindent The nine hyperplanes mentioned in Theorem 2 consist of matrices
of the form

\begin{equation} {\Pi}_1 = \left[ \begin{array}{cccc} B_{00} & B_{01} & \bullet & \bullet \\
B_{10} & B_{11} & \bullet & \bullet \\ \bullet & \bullet & \bullet & \bullet \\
\bullet & \bullet & \bullet & \bullet \end{array} \right] \hspace{6mm}
{\Pi}_2 = \left[ \begin{array}{cccc} B_{00} & B_{01} & \bullet & \bullet \\
\bullet & \bullet & \bullet & \bullet \\ B_{20} & B_{21} & \bullet & \bullet \\
\bullet & \bullet & \bullet & \bullet \end{array} \right] \end{equation}

\begin{equation} {\Pi}_3 = \left[ \begin{array}{cccc} B_{00} & B_{01} & \bullet & \bullet \\
\bullet & \bullet & \bullet & \bullet \\ \bullet & \bullet & \bullet & \bullet \\
B_{30} & B_{31} & \bullet & \bullet \end{array} \right] \hspace{6mm}
{\Pi}_4 = \left[ \begin{array}{cccc} B_{00} & \bullet & B_{02} & \bullet \\
B_{10} & \bullet & B_{12} & \bullet \\ \bullet & \bullet & \bullet & \bullet \\
\bullet & \bullet & \bullet & \bullet \end{array} \right]
\end{equation}

\begin{equation} {\Pi}_5 = \left[ \begin{array}{cccc} B_{00} & \bullet & B_{02} & \bullet \\
\bullet & \bullet & \bullet & \bullet \\ B_{20} & \bullet & B_{22} & \bullet \\
\bullet & \bullet & \bullet & \bullet \end{array} \right] \hspace{6mm}
{\Pi}_6 = \left[ \begin{array}{cccc} B_{00} & \bullet & B_{02} & \bullet \\
\bullet & \bullet & \bullet & \bullet \\ \bullet & \bullet & \bullet & \bullet \\
B_{30} & \bullet & B_{32} & \bullet \end{array} \right] \end{equation}

\begin{equation} {\Pi}_7 = \left[ \begin{array}{cccc} B_{00} & \bullet & \bullet & B_{03} \\
B_{10} & \bullet & \bullet & B_{13} \\ \bullet & \bullet & \bullet & \bullet \\
\bullet & \bullet & \bullet & \bullet \end{array} \right] \hspace{6mm}
{\Pi}_8 = \left[ \begin{array}{cccc} B_{00} & \bullet & \bullet & B_{03} \\
\bullet & \bullet & \bullet & \bullet \\ B_{20} & \bullet & \bullet & B_{23} \\
\bullet & \bullet & \bullet & \bullet \end{array} \right] \end{equation}

\begin{equation} {\Pi}_9 = \left[ \begin{array}{cccc} B_{00} & \bullet & \bullet & B_{03} \\
\bullet & \bullet & \bullet & \bullet \\ \bullet & \bullet & \bullet & \bullet \\
B_{30} & \bullet & \bullet & B_{33} \end{array} \right] \end{equation}

\noindent where the matrix elements that are explicitly written are assumed to
sum to one (and similarly for the remaining three blocks taken separately).

The normal vectors of these hyperplanes are the matrices

\begin{equation} n_1 = \frac{1}{4}\left[ \begin{array}{cccc} 1 & 1 & -1 & -1 \\
1 & 1 & -1 & -1 \\ -1 & -1 & 1 & 1 \\ -1 & -1 & 1 & 1 \end{array} \right]
\end{equation}

\noindent and so on.

\vspace{1cm}

{\bf Appendix B: Numerics}

\vspace{5mm}

\noindent {\bf I}. Average entropy in ${\cal B}_3$.

 To generate a random bistochastic matrix according to the
   flat measure on ${\cal B}_3\subset {\mathbb R}^4$, we have drawn
   at random a point $(x,y,z,t)$ in the $4$-dimensional hypercube.
   It determines a minor of a $N=3$ matrix, and the remaining
   five elements of $B_3$ may  be  determined by the
   unit sum conditions in eq. (\ref{1}). Condition
i is fulfilled if the sums in both rows and both columns of
the minor
    does not exceed unity, and the sum of all
 four elements is not smaller than one.
   If this was the case, the random matrix $B_3$ was accepted to the
ensemble of
   random bistochastic matrices. If additionally, the
   chain links condition (\ref{24}) were satisfied, the matrix was accepted
   to the ensemble of unistochastic matrices, generated with respect to the
   flat measure on ${\cal U}_3$. The mean entropies,  (\ref{entr3}),
    were computed by taking an  average over both ensembles
   consisting of $10^7$ random matrices, respectively.

\medskip

\noindent {\bf II}  Numerical verification, whether a given bistochastic matrix
$B$ is unistochastic.

   We have performed a random walk in the space of unitary matrices.
   Starting from an arbitrary random initial point $U_0$
  we computed $B_0=U_0\circ U^*_0$ and its distance
 to the analyzed matrix, $D_0=D(B_0,B)$, as defined in (\ref{dist}).
We were fixing a small parameter $\alpha \approx 0.1$,
generated a random Hermitian matrix $H$
distributed according to the Gaussian unitary ensemble   \cite{Mehta},
and found a unitary perturbation
$V=\exp(-\alpha H)$. The matrix $U_{n+1}=V U_n$
was accepted as a next point of the random trajectory,
if the  distance $D_{n+1}$ was smaller than the previous one, $D_n$.
If a certain number (say $100$) of random matrices $V$ did not
allow us to decrease the distance, we were reducing the angle $\alpha$ by
half, to start a finer search. A single run was stopped if
the distance $D$ was smaller then $\epsilon= 10^{-6}$ (numerical solution
found),
or $\alpha$ got smaller then a fixed cut off value (say
$\alpha_{\rm min}=10^{-4}$).
In the latter case, the entire procedure was repeated a hundred times,
starting from various unitary
random matrices $U_0$, generated according to the
Haar measure on $U(4)$ \cite{PZK98}. The smallest distance
$D_{\rm min}$ and the closest unistochastic matrix
$B_{min}=U_n\circ {\bar U_n}$ was recorded.

To check the accuracy of the algorithm
we constructed several random unistochastic matrices,
$B=U\circ {\bar U}$, and verified that random walk procedure
was giving their approximations with $D_{\rm min} < \epsilon$.

\vspace{1cm}

{\bf Appendix C: A system of equations}

\vspace{5mm}

In order to curtail a plethora of indices in Eqs. (\ref{37}-\ref{39}) and ease
the subsequent notation let us introduce shorthands: $\varphi_j=\phi_{j1}$,
$\psi_j=-\phi_{j2}$, $j=1,2,3$. With that the system rewrites as \vskip .2cm
\begin{eqnarray}
e^{i{\varphi}_{1}} + e^{i{\varphi}_{2}} + e^{i{\varphi}_{3}} = -L  \label{eqs1} \\
e^{i{\psi}_{1}} + e^{i{\psi}_{2}} + e^{i{\psi}_{3}} = -L \label{eqs2} \\
e^{i({\varphi}_{1}+ {\psi}_{1})} + e^{i({\varphi}_{2} + {\psi}_{2})} +
e^{i({\varphi}_{3} + {\psi}_{3})} =-L \label{eqs3},
\end{eqnarray}

\noindent We shall prove the following:

\

\noindent \underline{Lemma}: The system of equations (\ref{eqs1}-\ref{eqs3})

\noindent 1. has no real solutions for $L > 1$,

\noindent 2. for $0 < L < 1$ has the solution

\begin{equation} \begin{array}{lll} {\varphi}_{1} = 0 & {\varphi}_{2} = {\phi} &
{\varphi}_{3} = - {\phi} \\ {\psi}_{1} = -{\phi} & {\psi}_{2} = 0 & {\psi}_{3} =
{\phi} \end{array} \ , \hspace{6mm} \cos{\phi} = \frac{t-1}{t+1} =
- \frac{L+1}{2}  \ , \label{41a} \end{equation}

unique up to obvious permutations,

\noindent 3. for $L = 0,1$ has continuous families of solutions.

Indeed, each of the unimodal numbers $e^{i\varphi_k}$, $k=1,2,3$ is a root of:
\begin{eqnarray}
 P(\lambda) &=&(\lambda - e^{i\varphi_1})(\lambda  - e^{i\varphi_2})(\lambda  -
 e^{i\varphi_3})  \\ \nonumber &=&\lambda ^3 - (e^{i\varphi_1}  + e^{i\varphi_2}  +
 e^{i\varphi_3})\lambda ^2  + (e^{i(\varphi_1 + \varphi_2)}  + e^{i(\varphi_1 + \varphi_3)}  +
 e^{i(\varphi_2 + \varphi_3)})\lambda - e^{(i\varphi_1  + \varphi_2  + \varphi_3)} \\
 \nonumber &=&\lambda ^3 - (e^{i\varphi_1}  + e^{i\varphi_2}  + e^{i\varphi_3})\lambda
 ^2  + (e^{-i\varphi_3}  + e^{-i\varphi_2}  + e^{-i\varphi_1})e^{(i\varphi_1  + \varphi_2  +
 \varphi_3)}\lambda - e^{(i\varphi_1  + \varphi_2  + \varphi_3)} \\ \nonumber &=&\lambda
 ^3+\lambda^2L-\lambda Le^{i\Phi} -e^{i\Phi}=\lambda^2(\lambda+L)-(1+\lambda
 L)e^{i\Phi},
 \end{eqnarray}
where $\Phi=\varphi_1+\varphi_2+\varphi_3$, and we used (\ref{eqs1}) and the reality of
$L$. Thus each $\lambda = e^{i\varphi_k}$, $(k=1,2,3)$, fulfils:
\begin{equation}\label{condl}
\lambda^2(\lambda+L)=(1+\lambda L)e^{i\Phi}.
\end{equation}

Analogously, $\mu =e^{i\psi_k}$, $(k=1,2,3)$, fulfils
\begin{equation}\label{condm}
\mu^2(\mu+L)=(1+\mu L)e^{i\Psi},
\end{equation}
with $\Psi=\psi_1+\psi_2+\psi_3$.

Observe now, that if $\lambda=e^{i\varphi_k}$ and $\mu=e^{i\psi_k}$ are solutions
of (\ref{eqs1}-\ref{eqs3})
with the same number $k$, $(k=1,2,3)$ then, upon the same reasoning applied
to (\ref{eqs3}), $\lambda\mu$ fulfils
\begin{equation}\label{condlm}
\lambda^2\mu^2(\lambda\mu+L)=(1+\lambda\mu L)e^{i(\Phi+\Psi)} \ . 
\end{equation}

\noindent Multiplying (\ref{condl}) by (\ref{condm}) and finally by 
(\ref{condlm}) after exchanging its sides, we obtain, after division by
$\lambda^2\mu^2e^{i(\Phi+\Psi)}\ne 0$,
\begin{equation}\label{res1}
(L+\lambda)(L+\mu)(L\lambda\mu+1)=(L\lambda+1)(L\mu+1)(L+\lambda\mu),
\end{equation}
which, upon substitution $\lambda=e^{i\varphi_k}$, $\mu=e^{i\psi_k}$ and
putting everything on one side factorizes to
\begin{equation}\label{res2}
L(L-1)(e^{i\varphi_k}-1)(e^{i\psi_k}-1)(e^{i(\varphi_k+\psi_k)}-1)=0,
\end{equation}
(any computer symbolic manipulation program can be helpful in revealing
(\ref{res2}) from (\ref{res1})).

Hence, if $L\ne 0,1$, then for each pair $(\varphi_k,\psi_k)$, $k=1,2,3$,
either: a) one of the angles is zero or b) they are opposite. The latter case
can not occur for all three pairs since then $e^{i({\varphi}_{1}+ {\psi}_{1})}
+ e^{i({\varphi}_{2} + {\psi}_{2})} + e^{i({\varphi}_{3} +
{\psi}_{3})}=3\ne-L$, hence at least one of $\varphi_k$ or $\psi_k$ equals
zero. Up to unimportant permutations we can assume $\varphi_3=0$, but then,
since $e^{i{\varphi}_{1}} + e^{i{\varphi}_{2}} + e^{i{\varphi}_{3}} =-L\in
\mathbb{R}$, we immediately get $\varphi_1=-\varphi_2$. This determines also
all other angles (also up to some unimportant permutation) and we end up with
the solution announced in point 2. above as the only possibility, but such a
solution exists only if $L\le 1$.

To prove 3. observe that
\begin{enumerate}
\item for $L=0$,
\begin{eqnarray}\label{sol_0}
  &&\varphi_1=\varphi,\quad \varphi_2=\varphi+2\pi/3,\quad \varphi_3=\varphi+4\pi/3, \\
  &&\psi_1=\psi,\quad \psi_2=\psi+2\pi/3,\quad \psi_3=\psi+4\pi/3,
\end{eqnarray}
is a legitimate solution of (\ref{eqs1}-\ref{eqs3}) for arbitrary $\varphi$ and
$\psi$,
\item for $L=1$
\begin{eqnarray}\label{sol_1}
  &&\varphi_1=\varphi,\quad \varphi_2=\pi,\quad \varphi_3=\varphi+\pi \\
  &&\psi_1=-\varphi+\pi,\quad \psi_2=\pi,\quad \psi_3=-\varphi,
\end{eqnarray}
is a solution for an arbitrary $\varphi$.
\end{enumerate}


\begin{thebibliography}{99}

\bibitem{Sylvester} J. J. Sylvester, Phil. Mag. \underline{34} (1867) 461.

\bibitem{vanLint} J. H. van Lint, Math. Intelligencer \underline{4} (1982) 72.

\bibitem{Lande} A. Land\'e, \emph{From Dualism to Unity in Quantum Physics},
Cambridge U. P. (1960).

\bibitem{Rovelli} C. Rovelli, Int. J. of Theor. Phys. \underline{35} (1996) 1637.

\bibitem{Khrennikov} A. Khrennikov, J. Phys. \underline{A34} (2001) 1.

\bibitem{Hadamard} M. J. Hadamard, Bull. Sci. Math. \underline{17} (1893) 240.

\bibitem{Hedayat} A. Hedayat and W. D. Wallis, Ann. Stat. \underline{6} (1978) 1184.

\bibitem{Sloane} N. J. A. Sloane's homepage, www.research.att.com/\~{}njas/hadamard/
index.html.

\bibitem{Zeilinger} A. Zeilinger, M. Zukowski, M. A. Horne, H. J. Bernstein and D. M.
Greenberger, in J. Anandan and J. L. Safko (eds): \emph{Fundamental Aspects of Quantum
Theory}, World Scientific, Singapore (1994).

\bibitem{Stenholm} P. T\"orm\"a, S. Stenholm and I. Jex, 
Phys. Rev. \underline{A52} (1995) 4853.

\bibitem{Werner} R. F. Werner, J. Phys. \underline{A34} (2001) 7081.

\bibitem{Bjorck} G. Bj\"orck and B. Saffari, C. R. Acad. Sci. Paris, S\'er. I
\underline{320} (1995) 319.

\bibitem{Haagerup} U. Haagerup, in \emph{Operator Algebras and Quantum Field Theory,
Rome (1996)}, Internat. Press, Cambridge, MA (1997).

\bibitem{Petrescu} M. Petrescu, \emph{Existence of continuous families of complex
Hadamard matrices of certain prime dimensions and related results}, UCLA thesis,
Los Angeles (1997).

\bibitem{Dita} P. D\u{i}\c{t}a, arXiv: quant-ph/0212036.

\bibitem{Tanner} G. Tanner, J. Phys. \underline{A34} (2001) 8485.

\bibitem{Pakonski1} P. Pako\'nski, G. Tanner and K. \.Zyczkowski, J. Phys.
\underline{A34} (2001) 9303.

\bibitem{Pakonski2} P. Pako\'nski, G. Tanner and K. \.Zyczkowski, J. Stat.
Phys. \underline{111} (2003) 1331.

\bibitem{Karol} K. \.Zyczkowski, M. Ku\'s, W. S{\l}omczynski and H.-J. Sommers,
J. Phys. \underline{A36} (2003) 3425.

\bibitem{Jarlskog} C. Jarlskog and R. Stora, Phys. Lett. \underline{B208} (1988) 268.

\bibitem{Auberson} G. Auberson, A. Martin and G. Mennessier, Commun. Math. Phys.
\underline{140} (1991) 523.

\bibitem{Mennessier} G. Mennessier and J. Nyuts, J. Math. Phys. \underline{15}
(1974) 1525.

\bibitem{Tadej} W. Tadej et al., to appear.

\bibitem{Birkhoff} G. Birkhoff, Univ. Nac. Tucum\'an Rev. \underline{A5} (1946) 147.

\bibitem{Brualdi} R. A. Brualdi and P. M. Gibson, J. Comb. Theory \underline{A22}
(1977) 194.

\bibitem{WS} W. S{\l}omczy\'nski, Open Sys. Inf. Dyn. \underline{9} (2002) 201.

\bibitem{Beck} M. Beck and D. Pixton, arXiv: math.CO/0305322.

\bibitem{Poon} Y.-H. Au-Yeung and Y.-T. Poon, Lin. Alg. Appl. \underline{27} (1979) 69.

\bibitem{India} H. G. Gadiyar, K. M. S. Maini, R. Padma and
H. S. Sharatchandra, J. Phys. \underline{A36} (2003) L109.

\bibitem{Jones} K. R. W. Jones, J. Phys. \underline{A32} (1990) L1247.

\bibitem{Au-Yeung} Y.-H. Au-Yeung and C.-M. Cheng, Lin. Alg. Appl.
\underline{150} (1991) 243.

\bibitem{Mehta} M. L.   Mehta, {\sl Random Matrices}, II ed.,
New York: Academic (1991).

 \bibitem{PZK98} M.  Po{\'z}niak, K. {\.Z}yczkowski, and
M. Ku{\'s}, J. Phys. \underline{A31} (1998) 1059.

%\bibitem{Cencov} N. N. {\v C}encov, Statistical Decision Rules and Optimal
%Inference, Amer. Math. Soc., Providence 1982.

%\bibitem{Manin} Y. Manin

%\bibitem{Weinhold} F. Weinhold, J. Chem. Phys. \underline{63} (1975) 2479.

%\bibitem{Ihrig} P. Salamon, J. D. Nulton and E. Ihrig, J. Chem. Phys.
%\underline{80} (1984) 436.

%\bibitem{BTZ} M. Ba\~nados, C. Teitelboim and J. Zanelli,
%Phys. Rev. Lett. \underline{69} (1992) 1849.

%\bibitem{Gibbons} S. Ferrara, G. W. Gibbons and R. Kallosh, Nucl. Phys.
%\underline{B500} (1997) 75.

%\bibitem{BTZH} Ba\~nados, M. Henneaux, C. Teitelboim and
%J. Zanelli, Phys. Rev. \underline{D48} (1993) 1506.

%\bibitem{Davies} P. C. W. Davies, Proc. R. Soc. Lond. \underline{A353}
%(1977) 499.

%\bibitem{Hawking} S. W. Hawking and D. N. Page, Comm. Math. Phys.
%\underline{87} (1983) 577.

%\bibitem{Myers} A. Chamblin, R. Emparan, C. V. Johnson and R. C. Myers,
%Phys. Rev. \underline{D60} (1999) 104026.


\end{thebibliography}
\end{document}